\documentclass[12pt]{article}
\usepackage{theorem, amsfonts,amsmath,graphicx}
\usepackage{setspace}        

\nonstopmode
\newcommand\qed{\hfill$\sqcap\kern-8.0pt\hbox{$\sqcup$}$}

\title{Nonlinearity, correlation and the valuation of employee  options}

\author{M. R. Grasselli \\ 
McMaster University}

\begin{document}
\maketitle

\begin{abstract}
We propose a discrete time algorithm for the valuation of employee stock options based on exponential indifference prices and taking into account both the possibility of partial exercise of a fraction of the options and the use of a correlated traded asset to hedge part of their risk. We determine the optimal exercise policy under this conditions and present numerical results showing how both effects can significantly change the value of the option for an employee, as well as its cost for the issuing firm.
\end{abstract}


\section{Introduction}

Stock options given to employees as part of their compensation packages have received increasing amounts of academic attention due to recent public debates on regulations for their accounting status. Most of the debate concerns the lack of agreement on a consistent and tractable model for the valuation of such options. Typically, employee options consist of American calls issued by the company on its own stock, and carry several contractual features, such as a vesting period in which the option cannot be exercised, the forfeit of the option if the employee leaves the company during the life of the option, and trading restrictions on both the option itself and the underlying stock. 

In the absence of these contractual features, such options should be priced according to the standard paradigms of arbitrage and replication of derivatives in complete markets. For instance, assuming a geometric Brownian motion for the continuous time dynamics of the (non-dividend paying) underlying stock would lead to them being priced by the Black--Scholes formula, whereas a binomial model approximation for a discrete time dynamics (even with dividends) leads to the well-known Cox--Rubstein-Ross valuation algorithm. In the presence of the contractual features, it is clear that these models should fail to produce any sensible price for such options, since their basic hypotheses, such as the possibility of constructing a self-financing replicating portfolio, are violated. Our point of view is that the contractual features of employee options introduce risks that cannot be hedge in the financial market available to their holders, such as the risk of an employee leaving the company, or the market risk that cannot be hedged away by trading in the underlying stock. Therefore, employee options are financial derivatives in incomplete markets and should be priced as such. 

Not surprisingly, empirical research by Huddart and Lang (1996) and Carpenter (1998) confirms that the holders of such options behave differently than the holders of otherwise identical but unrestricted options. In general terms, they tend to exercise the options sooner and at a lower price for the underlying stock. Since this represents an opportunity loss when compared to the case of regular American call options, the corresponding Black--Scholes price (or that obtained in a binomial model) necessarily overestimates the value of employee options. This is in accordance with the common perception of executive, accountants, compensations consultants and some politicians, who reacted with fierce criticism to the proposal made in 1993 by the Financial Accounting Standard Board (FASB) on the way companies account for compensation packages. In its statement 123 of 1995,  FASB proposed that companies should deduct from corporate earnings the Black--Scholes value of employee options at grant-date. As reported in Hall and Murphy (2002) the general reaction from all interested parties was that ``the Black--Scholes price is too high".     

To address these criticisms, in 2004 FASB revised its statement 123, recommending the use of ``option-pricing models adjusted for the unique characteristics of those instruments". Specifically, it recommends that the (shorter) expected life of the options be estimated from historical data and used in place of the contractual life for valuation purposes. In addition, it recommends that the option value should be adjusted to account for the possibility of the employee leaving the company during the vesting period (and therefore forfeiting the option). In a recent paper, Hull and White (2004) observe that the FASB recommendations are not theoretically sound. They argue that the simple adjustment for the life of the option is unlikely to capture the ``true" value for the option since it does not take into account the employee exercise strategy. Instead, they propose to modify the standard binomial algorithm by explicitly assuming that employees exercise their options as soon as the underlying stock price reaches a threshold $M$, which should then be estimated from empirical data. After Hull and White (2004), other papers have appeared in the literature with analytic solutions for its  continuous time versions. For instance, Sircar and Xiong (2005) explore the infinite maturity limit within the framework of perpetual call options with an exercise barrier given by $M$. For finite maturities, Cvitanic, Wiener and Zapatero (2004) use 
the distributional properties of first passage times for a geometric Brownian motion to obtain explicit formulas for the value of employee options assuming that exercise occurs when the underlying stock hits a (time dependent) barrier $M_t=Me^{\alpha t}$. In all cases, the fact that employees might leave the company during the life of the option is model as an independent Poisson process and its effect on the value of the option is then calculated by techniques akin to the valuation of derivatives subject to default risk. 

Let us ignore for a moment the effect of employee exit rate and concentrate on the exercise barrier feature of these models. It rests on the assumption that {\em all} of the options in the compensation package should be exercised once the underlying stock price hits the barrier. Observe that this is a relic of the standard valuation paradigms for early exercise derivatives in complete markets, where prices are given by risk-neutral expectations and are therefore linear with respect to the size of the derivative. In particular, in a complete market, the decision to exercise one American call is identical to the decision to exercise any number of American calls. However, taking the point of view that employee options are derivatives in incomplete markets, we need to abandon the risk-neutral expectation paradigm and use the employee's risk preferences in order to find the exercise strategy. Since risk preferences, in whatever formulation one might adopt, are generally nonlinear, there is no justification in assuming that the exercise decision is independent of the number of options held by the employee. 

In fact, Rogers and Scheinckman (2003) use a utility based framework to show that the optimal exercise of American claims in incomplete markets depends explicitly on the number of units of the claim that are held at a given time and, moreover, consists of gradually exercising fractions of the holdings during the life of the options as opposed to exercising  them all at the same time.  For the particular case of employee options, this partial exercise effect was investigated further by Jain and Subramanian (2004), where the value that the options have for the employee is found by the certainty-equivalent approach previously applied to this problem by Lambert, Larker and Verrechia (1991) and Detemple and Sundaresan (1999). Their results unequivocally show that requiring all options to be exercised at once significantly underestimates both their value to the employee and their cost for the issuing company.  

In this paper, we extend this utility-based valuation of employee options to the case where trading on a correlated asset is allowed. From a practical point of view, this corresponds to the concrete situation where the employee is forbidden to trade on the underlying company stock, but is allowed to trade on an index with a known correlation with the company. As we will see, this can greatly increases the option value for the employee. Mathematically, this is the standard setup for a derivative written on a non-traded asset which can be partially hedge with a traded correlated asset. 

We propose to use an exponential indifference pricing framework to obtain the price of a derivative in this type of two-factor incomplete markets. For the case of European style claims, a discrete-time algorithm for indifference prices in these markets was studied by Musiela and Zariphopoulou (2004). When the claim allows for early exercise, the theoretical definition of its indifference price in continuous time was given in Oberman and Zariphopoulou (2003), who also analyze numerical schemes for viscosity solutions to the corresponding HJB type equations. In the limit of an infinite time horizon, an analytic solution for the case of a single American call option on a non-traded asset was obtained by Henderson (2005). In this paper, we describe a discrete time algorithm for the valuation and optimal exercise policy of {\em multiple} American call options with a {\em finite} time horizon in the presence of a 
{\em correlated} traded asset.  

In order to concentrate on the combined effects of nonlinearity and correlation on the valuation of employee options, we make the simplifying assumptions that the options have no vesting period and that the employee exit rate from the company is zero. In the last section we comment on how these two features can be naturally incorporated into our algorithm. The remaining of the paper is organized as follows. In section \ref{one-period} we explicitly calculate the indifference price, the exercise threshold and the optimal hedge strategy for American call options written on a non-traded asset in the presence of another correlated asset. In particular we show how the exercise threshold depends on the number of options held by the employee and how partial exercise at the initial time can arise as the optimal strategy.  In section \ref{two-period} we extended the analysis to a two-period model, where both partial and inter-temporal exercise are allowed and show how to obtain the indifference price for multiple American call options locally in the tree. The setup is then generalized to a multi-period model in \ref{multi-period}, where the option values and optimal exercise policy are calculated along a rectangular grid of time steps and stock values. Using this multi-period algorithm, numerical results for the optimal exercise policy, employee option value and the effective cost for the firm are presented in section \ref{numerical}.  

\section{The one-period model: volume dependent exercise}
\label{one-period}

Consider a one-period model for a market consisting of two risky assets and a riskless cash account. The cash account is assumed to have a constant {\em discounted} value normalized to 1, whereas the risky assets have {\em discounted} values given by the two-dimensional random variable 
\begin{equation}
(S_T,Y_T)=\left\{\begin{array}{l}
(uS_0,hY_0) \quad \mbox{ with probability } p_1, \\
(uS_0,\ell Y_0) \quad \mbox{  with probability } p_2, \\
(dS_0,h Y_0) \quad \mbox{ with probability } p_3, \\
(dS_0,\ell Y_0) \quad \mbox{  with probability } p_4, \end{array}\right.
\end{equation}
where $0<d<1<u$ and $0<\ell<1 <h$, for positive initial values $S_0,Y_0$ and historical probabilities $p_1,p_2,p_3,p_4$.

We assume that an investor in this market can keep and amount $\xi$ in the cash account, while holding $H$ units of the asset $S$, but is not allowed to invest in the correlated asset $Y$. His wealth then evolves from an initial value $x=\xi+HS_0$ to a terminal discounted value
\begin{equation}
X^x_T=\xi+HS_T=x+H(S_T-S_0).
\end{equation}

We now introduce a contingent claim with maturity $T$ whose payoff $C_T \in {\cal F}$ is written as a function $C(Y_T)$ of the non-traded asset, therefore rendering the market incomplete.  Given an exponential utility function $U(x)=-e^{-\gamma x}$, an investor who {\em buys} this claim for a price $\pi$ will then try to solve the optimal portfolio problem
\begin{equation}
u^C(x-\pi)=\sup_{H} E[U(X^{x-\pi}_T+C_T)]=-e^{-\gamma (x-\pi)} \inf_{H} E\left[e^{-\gamma H(S_T-S_0)-\gamma C(Y_T)}\right].
\label{optimalhedge}
\end{equation}

The {\em indifference price} for this claim is defined to be the amount $\pi^C$ that makes the investor indifferent between buying the claim and then optimally hedge it according to \eqref{optimalhedge} or foregoing the deal and just optimally manage his wealth $x$. In other words, it is the unique solution to the equation
\begin{equation}
u^0(x)=u^C(x-\pi),
\label{ind_id}
\end{equation}
where $u^0$ is defined by \eqref{optimalhedge} for the degenerate case $C\equiv 0$.

Denoting the possible pay-offs for the derivative $C$ at the terminal time by $C_h=C(hY_0)$ and $C_\ell=C(\ell Y_0)$, it is straightforward to obtain that the number of shares $H^C$ in the optimal hedging portfolio solving \eqref{optimalhedge} must satisfy
\begin{equation}
\label{Hoptimal}
e^{-\gamma (u-d)H^CS_0}=\frac{p_3e^{-\gamma C_h}+p_4e^{-\gamma C_\ell}}{p_1e^{-\gamma C_h}+p_2e^{-\gamma C_\ell}}\left(\frac{1-d}{u-1}\right).
\end{equation}
For the case $C\equiv 0$, this gives 
\begin{equation}
\label{Hmerton}
e^{-\gamma (u-d)H^0S_0}=\frac{p_3+p_4}{p_1+p_2}\left(\frac{1-d}{u-1}\right),
\end{equation}
corresponding to Merton's optimal investment portfolio. Therefore, the {\em excess hedge} that is necessary to account for 
the claim $C$ is given by $H^C-H^0$. 

Inserting these expressions into the value functions $u^0$ and $u^C$, we obtain the following expression for the indifference price
\begin{equation}
\label{european}
\pi^C=g(C_h,C_\ell)
\end{equation}
where, for fixed parameters $(u,d,p_1,p_2,p_3,p_4)$ the function $g:\mathbb{R}\times\mathbb{R}\rightarrow \mathbb{R}$ is given by 
\begin{equation}
g(x_1,x_2)=\frac{q}{\gamma}\log \left(\frac{p_1+p_2}{p_1e^{-\gamma x_1}+p_2e^{-\gamma x_2}}\right)+
\frac{1-q}{\gamma}\log\left(\frac{p_3+p_4}{p_3e^{-\gamma x_1}+p_4e^{-\gamma x_2}}\right),
\label{gfunction}
\end{equation}
with
\[q=\frac{1-d}{u-d}.\]
Observe that for either perfectly correlated or perfectly anti-correlated assets, corresponding respectively to $p_2=p_3=0$ or $p_1=p_4=0$, the market is reduced to a complete one and the expression above correctly gives the risk-neutral price for $B$, that is
\begin{equation}
\pi^C=q C(hY_0) + (1-q) C(\ell Y_0) =: E^Q[C(Y_T)],
\label{complete}
\end{equation}
where $Q$ is the unique martingale measure for this complete market. There is another way to reduce \eqref{european} to an expectation functional, namely  
\begin{equation}
\lim_{\gamma \rightarrow 0} \pi^C=(\tilde q_1+\tilde q_3) C(hY_0)+(\tilde q_2+\tilde q_4)C(\ell Y_0)=:E^{\tilde Q}[C(Y_T)],
\label{minimal}
\end{equation}
where 
\begin{eqnarray}
\tilde q_1&=&q\frac{p_1}{p_1+p_2}, \quad \qquad \quad \tilde q_2=q\frac{p_2}{p1+p_2} \nonumber \\
\tilde q_3&=&(1-q)\frac{p_3}{p_1+p_4}, \quad \tilde q_4=(1-q)\frac{p_4}{p_3+p_4}
\end{eqnarray} 
defines the so-called {\em minimal martingale measure} for this market. For a general incomplete market and non-zero risk aversion, the indifference price does not correspond to an expectation with respect to a martingale measure. However, Musiela and Zariphopoulou (2004) show how to rewrite \eqref{european} in terms of a nonlinear pricing operator making use of 
$\tilde Q$. They provide further economic insight for this pricing procedure by separating its {\em actuarial} from its {\em arbitrage-free} features. For the purposes of our paper, however, expression \eqref{european} will suffice, and the reader is referred to their paper for additional interpretations of their pricing functional.

To complete the discussion on the properties of expression \eqref{european}, we can differentiate it with respect to $\gamma$ to verify that, provided  $C(Y_T)\geq 0$, the indifference price is a {\em decreasing} function of the risk aversion parameter. In particular
\begin{equation}
\lim_{\gamma \rightarrow \infty} \pi^C=C(\ell Y_0),
\end{equation}
corresponding to the smallest possible arbitrage-free price for the claim $C$, also known as the {\em super-replication price}.  By taking the second derivative of the indifference price with respect to 
$\gamma$, we find that $\pi^C$ is a {\em concave} function of the risk aversion. Moreover, we can use that scaling property
\[\pi^{aC}(\gamma)=a\pi^C(a \gamma)\]  
to conclude that the map 
\begin{equation}
a\rightarrow \frac{1}{a}\pi^{aC}
\end{equation} is a decreasing function of $a$. In other words, the {\em per unit} indifference 
price of multiple units of the same option {\em decreases} as the number of options increases.   

Suppose now that the contract $C$ allows for early exercise, so that its holder either waits until the terminal time $T$ or exercises it immediately at $t=0$.
It is clear that early exercise will occur whenever 
\[C(Y_0)\geq \pi,\]
where $\pi$ is the indifference price of an European claim with an equivalent pay-off at time $T$.
For example, it is straightforward to verify that an American call option on $Y$ with strike price $K$ will be exercised at time $t=0$ if and only if $Y_0$ exceeds the solution to the 
equation
\begin{equation}
(Y-K)^+=g((hY-K)^+,(\ell Y-K)^+),
\end{equation}
where $g:\mathbb{R}\times\mathbb{R}\rightarrow \mathbb{R}$ is given by \eqref{gfunction}.

As we have just seen, the indifference price for the buyer of a contingent claim scales sub-linearly with the size of the claim. 
This is a direct consequence of the fact that a risk averse investor who is prepared to pay $\pi$ dollars for a risky contract will not be convinced to pay a thousand times as much for a thousand units of the same contract. As a result, the early exercise threshold for an American call option obtained above is different (and higher) than the exercise threshold for a contract consisting of $A$ units of identical Americal calls. Explicitly, it it is obtained as the solution to 
\begin{equation}
A(Y-K)^+=g(A(hY-K)^+,A(\ell Y-K)^+)
\end{equation}

In figure 1, we show this threshold as a function of the number of options for $K=2$ and three different levels of correlation, with risk aversion parameter set to $\gamma=1$. The straight line corresponds to perfectly correlated assets which, as we have seen, reproduces risk-neutral valuation. The steeper line corresponds to uncorrelated assets, while the intermediate line corresponds to a correlation $\rho= 0.9554.$

Therefore, if compelled to exercise all his options at once, an investor holding multiples calls with the same strike price will do so at much reduced thresholds. If partial exercise is allowed, then the optimal number of options to be early-exercised  is
\begin{equation}
a_0=\arg \max_{0\leq a\leq A}\left[a(Y_0-K)^++\pi^{(A-a)}\right]
\label{optimal_number}
\end{equation}
where $\pi^{(A-a)}$ is the indifference price of an European option with pay-off 
\[C^{(A-a)}_T=(A-a)(Y_T-K)^+.\] 
If the indifference price were a linear operator, it is clear that the maximum for this expression would occur either at $a=0$ (no exercise) or at $a=A$ (full exercise), depending only on the relative slopes of two straight lines. Its nonlinearity, however, allows for the possibility of a maximum occurring at partial exercise.
This is depicted in a concrete example in figure 2 for $K=2$, $\gamma=0.2$ and correlation $\rho=0.2$. Observe that in this case the optimal policy consists of immediately exercising only 8 out of 10 American options, while keeping 2 until maturity.

Having determined $a_0$, we conclude that the value at time zero of $A$ units an American call option is
\[C^{(A)}_0=a_0(Y_0-K)^++\pi^{(A-a_0)}.\]

\section{Inter-temporal exercise in a two-period model}
\label{two-period}

An immediate generalization is to consider a model where the time interval $[0,T]$ is divided into a partition 
$0=t_0\leq t_1 \leq t_2=T$. Starting at an initial value $(S_0,Y_0)$, the two-dimensional stochastic process $(S_n,Y_n)$, for $n=1,2$, is then defined by the recursive relation 
\begin{equation}
(S_n,Y_n)=\left\{\begin{array}{l}
(uS_{n-1},hY_{n-1}) \quad \mbox{ with probability } p_1, \\
(uS_{n-1},\ell Y_{n-1}) \quad \mbox{  with probability } p_2, \\
(dS_{n-1},h Y_{n-1}) \quad \mbox{ with probability } p_3, \\
(dS_{n-1},\ell Y_{n-1}) \quad \mbox{  with probability } p_4, \end{array}\right.
\label{one-period-prob}
\end{equation}
where, as before, $0<d<1<u$ and $0<\ell<1 <h$, and $p_1,p_2,p_3,p_4$ are historical probabilities for the one-period model.

We now want to obtain the optimal exercise policy for multiple American calls knowing that they can be partially exercised at times $t_0$ and $t_1$. For each realization of the asset values, we seek to maximize an expression of the form \eqref{optimal_number}, where some of the options are exercised and some are kept alive. Once a certain number of options is exercised, their value equal their total pay-off, while the remaining portion is deem to have a value equal to the combined (non-linear) indifference price of equivalent number of European contracts. 

Observe that both the exercise value, given by a multiple of the pay-off $C(Y)$, and the continuation value, given by an indifference price $\pi$ of the form \eqref{european}, for options written on the non-traded asset depend only on the realized values of $Y_t$. The dependence on traded asset $S_t$ is only implicitly given in expressions such as \eqref{european} through the dynamic parameters $u,d$ and the corresponding historical probabilities. Therefore, we can reduce the dimensionality of our tree and concentrate exclusively on the time evolution of the non-trade asset $Y_t$.  

Let us denote by $h$ the node where the non-traded asset has value $hY_0$. If we find ourselves holding $A$ options at this node, then the optimal number of options to be exercised immediately is   
\begin{equation}
a_{h}=\arg \max_{0\leq a \leq A}\left[a(hY_0-K)^++\pi_{h}^{(A-a)}\right],
\label{optimal_number_high}
\end{equation}
where $\pi_{h}^{(A-a)}$ denotes the indifference of an European claim to starting at the node $h$ and maturing at time $T$ with pay-offs
\[C^{(A-a)}_{hh}=(A-a)(hhY_0-K)^+\]
with probability $(p_1+p_3)$ and 
\[C^{(A-a)}_{h\ell}=(A-a)(h\ell Y_0-K)^+\]
with probability $(p_2+p_4)$, where we have used $hh$ and $h\ell$ to denoted, respectively,  the nodes where the non-traded asset has values $hhY_0$ and $h\ell Y_0$.
Similarly to 
\eqref{european}, such indifference price is explicitly given by
\begin{equation}
\pi^{(A-a)}_{h}=g(C^{(A-a)}_{hh},C^{(A-a)}_{h\ell})
\end{equation}
In the same vein, the optimal number of options to be exercised at the node $\ell $, where the non-trade asset has value 
$\ell Y_0$, is
\begin{equation}
a_{\ell}=\arg \max_{0 \leq a \leq A}\left[a(\ell Y_0-K)^++\pi_{\ell}^{(A-a)}\right],
\label{optimal_number_low}
\end{equation}
where 
\begin{equation}
\pi^{(A-a)}_{\ell}=g(C^{(A-a)}_{h\ell},C^{(A-a)}_{\ell\ell})
\end{equation}

Therefore, at the intermediate time $t_1$, the total value of $A$ options at the node $h$ is
\begin{equation}
C^{(A)}_{h}:= \left[a_{h}(hY_0-K)^++\pi_{h1}^{(A-a_{h})}\right],
\end{equation}
while the total value of $A$ options at the node $\ell$ is
\begin{equation}
C^{(A)}_{\ell}:= \frac{1}{A}\left[a_{\ell}(\ell Y_0-K)^++\pi_{\ell 1}^{(A-a_{\ell })}\right].
\end{equation}

The final step in our valuation algorithm in the two-period model consists in finding the optimal number of options to be exercised at time $t_0=0$, that is,
\begin{equation}
a_0=\arg \max_{0\leq a\leq A}\left[a(Y_0-K)^++\pi_0^{(A-a)}\right],
\label{optimal_number_initial}
\end{equation}
where $\pi_0^{(A-a)}$ is the indifference price for of an European claim starting at time $t_0$ and maturing at time $t_1$ with possible pay-offs equal to $C^{(A-a)}_{h1}$, with probability $(p_1+p_3)$, and $C^{(A-a)}_{\ell 1}$, with probability $(p_2+p_4)$. Appealing to our previous calculations, this is given by 
\begin{equation}
\pi^{(A-a)}_{h}=g(C^{(A-a)}_{h},C^{(A-a)}_{\ell})
\end{equation}
Therefore, according to this valuation scheme, the value at time zero of $A$ units of an American call option on the non-traded asset is 
\begin{equation}
C^{(A)}_0:= \left[a_0(Y_0-K)^++\pi_0^{(A-a_0)}\right].
\end{equation}

\section{The multi-period model}
\label{multi-period}

In order to approximate a continuous time dynamics for the assets, we can consider a multi-period model  where the time interval $T$ is divided into $N$ subintervals with equal time steps $\Delta t = T/N$, with the dynamics for the asset prices on each subinterval given by \eqref{one-period-prob} for $n=1,\ldots,N$.  In the limit of small $\Delta t$, we want the discrete-time evolution to reproduce the dynamics of a two-factor Markovian market 
\begin{eqnarray}
dS_t &=& (\mu-r)S_tdt +\sigma S_t dW \\
dY_t &=& (\alpha-r-\delta)Y_tdt +\beta Y_t (\rho dW + \sqrt{1-\rho^2}) dZ,
\end{eqnarray}
where $\mu,\alpha,\sigma,\beta,r,\delta$ are constants and $(W,Z)$ are standard independent $P$--Brownian motions. We then choose the dynamic parameters $u,d,h,\ell$ and the one-period probabilities $p_i$, $i=1,2,3,4$, in order to match the distributional properties of the continuos time processes $S_t$ and $Y_t$. For instance, to obtain the correct variance, we can take 
\begin{eqnarray}
u &=& e^{\sigma\sqrt{\Delta t}}, \qquad h=e^{\beta\sqrt{\Delta t}} \\
d &=& e^{-\sigma\sqrt{\Delta t}}, \qquad \ell=e^{-\beta\sqrt{\Delta t}},
\end{eqnarray}
whereas for the correct mean we need
\begin{eqnarray}
p_1+p_2&=&\frac{e^{(\mu-r)\Delta t}-d}{u-d} \label{system1}\\
p_1+p_3&=&\frac{e^{(\alpha-r-\delta)\Delta t}-\ell}{h-\ell} \label{system2}
\end{eqnarray}
Finally, to obtain the desired correlation between $S_t$ and $Y_t$ we take
\begin{equation}
\rho b \sigma \Delta t = (u-d)(h-\ell)[p_1p_4-p_2p_3]. \label{system3}
\end{equation}
The solution to \eqref{system1},\eqref{system2} and \eqref{system3}, supplemented by the condition 
\begin{equation}
p_1+p_2+p_3+p_4=1
\end{equation}
uniquely determines the historical probabilities $p_i$.

For these parameters, the possible outcomes for the non-traded asset over $N$ time steps are 
\begin{equation}
Y^{(i)}=h^{N+1-i}Y_0, \qquad i=1,\ldots, 2N+1,
\label{array}
\end{equation}
ranging from $(h^NY_0)$  to $(\ell^NY_0=h^{-N}Y_0)$, respectively the highest and lowest achievable values 
starting from the middle point $Y_0$ with the multiplicative parameter $h=\ell^{-1}>1$. A realization for the discrete-time process $Y_n$ following the dynamics \eqref{one-period-prob} can then be thought as a path over a $(2N+1)\times N$ rectangular grid having 
the values in \eqref{array} as its columns. 

Our goal is to obtain an optimal exercise policy for the holder of $A$ units of an American call option with pay-off 
$(Y_n-K)^+$ over the time interval $[0,T]$, given that partial exercise of an integer number of these options is allowed at any intermediate time. We do this by determining the optimal number of options to be exercised at each node
$(in)$ of our rectangular grid, where the first index corresponds to the value $Y^{(i)}$ for the non-traded asset at the node, while the second index indicates the time $t_n$, $n=1,\ldots,N$. This will also determine the per unit value of the options on the entire grid, in particular at time zero, when they are awarded to the employee. 

At maturity, all options should be exercised whenever the terminal value $Y_T$ exceeds the strike price $K$. Therefore, using the values for the non-traded asset given in \eqref{array}, the optimal number of options to be exercised at time $t_N=T$ is given by
\begin{equation}
a_{iN}=\left\{\begin{array}{l}
A \qquad \mbox{if } Y^{(i)}\geq K, \\
0 \qquad \mbox{    otherwise},\end{array}\right. 
\label{terminal}
\end{equation}
for $i=1,\ldots,2N+1$. Accordingly, the value of $A$ options at the last column of the rectangular grid is
\begin{equation}
C^{(A)}_{iN}=A(Y^{(i)}-K)^+ , \qquad i=1,\ldots,2N+1.
\end{equation}
Next we specify the boundary conditions at the top and bottom of the grid (corresponding respectively to the highest and lowest values for $Y$) by 
\begin{eqnarray}
a_{1n}&=& A, \\
a_{2N+1,n}&=& 0,
\end{eqnarray}
for $n=0,\ldots,N$.

We can then obtain the optimal number of options to be exercised at intermediate times through the recursive relation
\begin{equation}
a_{in}=\arg \max_{0 \leq a \leq A, a\in\mathbb{Z}_+}\left[a(Y^{(i)}-K)^++\pi_{in}^{(A-a)}\right],  \quad \begin{array}{l}
n=N-1,\ldots, 0 \\ i=2,\ldots, 2N \end{array}
\end{equation}
where $\pi_{in}^{(A-a)}$ is the indifference price an European claim starting at the node $(in)$ with possible  
pay-offs $C^{(A-a)}_{i+1,n+1}$ and $C^{(A-a)}_{i-1,n+1}$, with probabilities $(p_1+p_3)$ and  $(p_2+p_4)$, respectively. From the previous sections, we know that this indifference price can be calculated explicitly as
\begin{equation}
\pi^{(A-a)}_{in}=g(C^{(A-a)}_{i+1,n+1},C^{(A-a)}_{i-1,n+1})
\end{equation}
Having calculated $a_{in}$, the value at the node $(in)$ of $A$ units of the option is given by
\begin{equation}
C^{(A)}_{in}=\left[a^*_{in}(Y^{(i)}-K)^++\pi_{in}^{(A-a_{in})}\right].
\end{equation}

\section{Numerical results}
\label{numerical}

\subsection{Exercise surface}

Let us first investigate the consequences of our model for the optimal exercise policy followed by an employee. 
With the notation of the previous section, let us define the {\em critical surface} as the number $(A-a_{in})$ of unexercised options that the employee should hold at time $t_n$ when the value of the non-traded asset is equal to $Y^{(i)}$ and the total number of options equal $A$. The optimal exercise policy then consists of exercising as many options as necessary to place the employee below this surface at any pair $(t_n,Y^{(i)})$. 

Consider the fixed parameters 
\begin{eqnarray*}
\mu &=& 0.12, \qquad \sigma=0.2, \qquad S_0=1.2, \qquad r=00.7\\
\alpha&=&0.15, \qquad \beta=0.3, \qquad Y_0=1.0, \qquad T=5\\
 N&=&500, \qquad A = 10, \qquad  K=1.0.
\end{eqnarray*}
Figure 3 depicts this critical surface for four different combinations of the remaining parameters. Our base case is shown in the top-left corner, for which the dividend rate is
$\delta=0.075$, the risk aversion is $\gamma=0.125$ and the correlation is $\rho=-0.5$. As expected, the optimal number of unexercised  options to be held at any given time decreases as the value of the underlying increases. However, this decrease does not occur as a sudden drop from $A$ to zero, which would characterize the exercise of all the options at once. Instead, this critical surface indicates that it is optimal for the employee to gradually unwind his holding. This preference for partial exercised is a direct consequence of the nonlinearity of risk preferences with respect to the number of outstanding options. It cannot be reproduced by any model which assumes that the value of the employee compensation package is linear in the number of options, such as the one proposed by Hull and White \cite{HulWhi04}. As we will see later, the possibility of partial exercise can significantly increase both the value that the options have for the employee and the cost that they represent for the firm.  

The preference for partial exercise is accentuated in the case shown in the top-right corner of figure 3, which differs from the base case by having $\delta=0$. In the absence of high dividends, it is optimal for the employee to wait for higher values of the underlying before exercising his options, and the exercise policy is even more gradual, corresponding to a less steep surface. 

Within our framework, a scenario leading to an optimal exercise policy resembling the constrain of exercising all options at once can be produce in two distinct ways. The first one is to impose a very high risk aversion, such as in the case depicted in the bottom-left corner of figure 3, which differs from the base case by having a high risk-aversion parameter $\gamma=10$. In this scenario, the valuation is still highly nonlinear, but the extremely low risk tolerance forces the employee to quickly exercise all his options as soon as the price of the underlying asset crosses a threshold. 

Alternatively, a critical surface dropping from $A$ to zero arises in the cases where the indifference price is given by a linear operator. According to \eqref{complete} and \eqref{minimal}, this occurs either for a market with perfectly correlated assets or for an employee with zero risk aversion. In the bottom-right corner of figure 3 we present the case of a  high-correlation $\rho=0.95$ between the traded and the non-traded assets, so that the market is close to completeness, therefore eliminating the nonlinearities inherent to indifference pricing. As a result, the exercise decision is independent of the number of options and the critical surface also exhibits a sudden drop from $A$ to zero.

\subsection{Option value for the employee}

We now consider how the option value predicted by our model varies with respect to the underlying parameters. In what follows, we want to simultaneously draw attention to the effects of nonlinearity of the risk preferences and correlation between the non-traded underlying stock $Y_t$ and the traded asset $S_t$. 

The nonlinear effects will favor partial exercise during the life of the options, resulting in an increased option value when compared with a model enforcing exercise of all options when the underlying reaches a volume independent threshold. For comparison, we calculate the option value in this {\em constrained} model 
by restricting the maximization at each node of our discrete grid to either $a=0$ or $a=A$ exercised options.  

As for correlation, the possibility of using a traded asset $S_t$ to partially hedge the risks carried by the claim $C(Y_t)$ using the excess hedge $(H^C-H^0)$ from \eqref{Hoptimal} and \eqref{Hmerton} has the effect of increasing the employee option value when compared with the case where no correlated asset is available for trading. Accordingly, we compare the option value predicted by our model for different levels of correlation, including the case $\rho=0$, corresponding to no excess hedge from the trading asset. For comparison with what can be obtained at high levels of correlation, we also compute the Black-Scholes price of the options for each combinations of the parameters, given that in the no-dividend case it corresponds to the correct complete market value of an American call option. In all of the subsequent pictures, the term option value refers to their {\em per unit value}, obtained by dividing the value of the entire compensation package by the number of options in it.

In figure 4, we present the dependence of the option value with time to maturity. The fixed parameters for these two plots are
\begin{eqnarray}
\mu &=& 0.09, \qquad \sigma=0.40, \qquad S_0=1.2, \qquad r=0.06 \nonumber \\
\alpha&=&0.08, \qquad \beta=0.45, \qquad Y_0=1.0, \qquad \delta=0 \nonumber \\
N&=&100, \qquad A= 10, \qquad  \quad  K=1.0. 
\label{par1}
\end{eqnarray} 
In both plots we observe that the option value increases with time to maturity, as expected for call options with early exercise. As time goes by, however, the employee option value increasingly deviates from the complete market Black--Scholes price, indicating that risk averse employees in face of hedging restrictions tend to exercise their options earlier than the corresponding option holder for whom a perfect hedge is available. Moreover, this early exercise effect is increased when the risk aversion is higher, as can be verified by comparing the two plots. The second relevant aspect of this picture is the confirmation that higher correlation leads to higher option values. Finally, we observe that the constrained model, even for a very high correlation 
$\rho=0.9$ produces lower option values than those obtained when partial exercise is allowed, especially for longer maturities. This is observe more acutely for a higher risk aversion, where the constrained model under-perfomrs its partial exercise counterpart even for very short maturities. 

Moving on to figure 5, we can analyze how the option value depends on the risk aversion parameter. The fixed parameters are the same as in \eqref{par1}. In accordance with the results derived for the one-period model, the option value is a decreasing concave function of risk aversion.  As before, the longer the maturity the larger the deviance from the Black--Scholes price, itself independent of risk preferences. Once more, higher correlations lead to higher option values. We also observe that while for short maturities and low risk aversion the option values generated by the constrained model are comparable to the those with partial exercise (albeit with a much higher  correlation) they quickly drop below all other partial exercise scenarios as soon as the risk aversion begins to increase.  

In figure 6 we observe the behavior of employee option values with changes in volatility of the underlying stock. We use the same fixed parameters as in \eqref{par1}, except that both $\sigma$ and $\beta$ are replaced by the same varying volatility. 
In general terms, we observe two conflicting trends: on the one hand, the employee's risk aversion favors early exercise (and reduced option value) as the risk (as measured by volatility) of the position increases; on the other hand, the convexity of the option pay-off leads to higher option values (and later exercise) for higher volatilities. This is clearly reflected in the four scenarios of figure \eqref{vol}. In the top-left corner, low risk aversion and short maturity makes the convexity of the option the predominant effect, and the employee option value increases with volatility, mimicking the behavior of the Black-Scholes price. In the bottom-left corner, an increased risk aversion parameter, keeping the same short maturity, produces comparable effects of risk aversion and convexity, resulting in practically constant option values with respect to volatility (expect for very high correlations). Moving to the bottom left-corner we see how, at longer maturities, the effect of risk aversion outweighs that of convexity and option values become a decreasing function of volatility. More interestingly, at the top-right corner we observe a combination of parameters for which the option value is not at all monotone with volatility: for some correlation it is a decreasing function, while for others it is decreasing up to a threshold, after which it starts increasing with volatility. This non-monotome behavior was also predicted by the model of Jain and Subramanian (2004), who pointed that for realistic parameters the volatility threshold would be too high to be observed. We are in agreement with these remarks, since their model corresponds to 
$\rho=0$ in our formulation. Our results show, however, that when trading in a correlated asset is allowed, the non-monotne behavior of option values can arise at realistic levels of volatility.        

In our last figure for this section (figureÊ7), we depict the dependence of employee option values with the correlation between the traded asset $S_t$ and the non-traded underlying $Y_t$. As expected from our previous analysis, option values increase as the possibility of partially hedging the risk with a correlated asset increases. Since a perfect hedge can be done with either a perfectly correlated or a perfectly anti-correlated asset, the option value increases symmetrically as the correlation moves away from zero, tending towards the Black--Scholes price at the limits $\rho\rightarrow \pm 1$. Once more, higher risk aversion and longer maturities produce larger deviance from the complete market value. The value obtained with the constrain that all 
options should be exercise at the same time presents the same behavior with respect to correlation. As always, for the same level of risk aversion, it consistently underperforms the model where partial exercise is allowed.  

\subsection{Cost for the firm}

While trading and hedging restrictions are imposed on the employee, we assume that the company granting the options in a compensation package is well diversified and can freely hedge them in the market. It then follows that, given an exercise policy for the employee, the total cost for the firm is obtained as the {\em discounted risk-neutral expected value} of the option pay-off at exercise date. In order to compute this cost, we first run our algorithm in order to determine the optimal exercise policy 
$a_{in}$ on a discrete grid for time and stock values. We then use Monte Carlo simulations of the {\em discounted risk-neutral} dynamics 
\begin{equation}
dY_t = -\delta Y_t dt + \beta Y_t dW^Q
\end{equation}
For each simulated value $Y_{in}$, we find the closest grid value $Y^{(i)}$ and determine the optimal number of options to be exercised at the corresponding grid node. We then add the discounted pay-offs of the exercised options over time and average over all the Monte Carlo paths. The results are shown in figure 8 using the following fixed parameters:
\begin{eqnarray}
\mu &=& 0.09, \qquad \sigma=0.40, \qquad S_0=1.2, \qquad r=0.06, \nonumber \\
\alpha&=&0.08, \qquad \beta=0.45, \qquad Y_0=1.0, \qquad \delta=0, \nonumber \\
N&=&100, \qquad A= 10, \qquad  \quad  K=1.0, \qquad T=5 \nonumber \\
\rho&=&0.6, \quad \quad T=5, \quad \qquad \gamma=0.5
\label{par2}
\end{eqnarray} 
We find that the effective cost for the firms lies between the value of the option from the point of view of the undiversified employee and the corresponding Black-Scholes price for the option. Moreover, its behavior with respect to variations in the model parameters is the same as the behavior of the employee option value, namely it increases with time to maturity, decreases with risk aversion, increases as the correlation with the traded asset moves away from zero and, for the choice of $T=5$ and $\gamma=0.5$, increases with volatility. 

\section{Conclusions}

We have formulated a discrete time algorithm for the valuation of employee options in the presence of a correlated traded asset and taking into account the possibility of partial exercise during the life of the options. The results obtained using the algorithm show that the nonlinearity of risk preferences leads to exercise strategies that depend on the number of outstanding option in a compensation package. Ignoring this dependence and imposing that all options should be exercised at the same time produces  systematic underestimates for the value that these options have to employee. Similarly, the possibility of hedging using a correlated asset greatly enhances the value of the options. Consequently, an issuing firm that neglects these two effects, nonlinearity and correlation, will estimate its costs much lower than they effectively are. 

We have shown how the valuation changes with the underlying model parameters. These results agree with the qualitative behavior predicted by the theoretical model of Jain and Subramanian (2004) in the absence of a correlated asset. They also agree with the empirical evidence in Huddart and Lang (1996) in the sense that employee options have much lower values than their corresponding Black--Scholes price, with increased discrepancy for longer maturity times and higher volatilities. Moreover, as expected intuitively, the options values decreases with the employee risk aversion and increase as the correlation with the traded asset increases. Finally, in accordance with Hall and Murphy (2002), we find that the effective cost for the issuing firm, although lower than the Black--Scholes price for the options, is higher than the value that the options have to employees. 

The extension of the model to allow for a vesting period should be immediate. Namely, we set the initial time in our algorithm to be $t_0=T_v$, where $T_v$ is the end of the vesting period, and calculate the value for the employee options for different levels of the underlying assets at time $T_v$. These values are then used as the pay-offs of an European style derivative with maturity $T_v$. We then use the same type of binomial tree valuation from the granting date to $T_v$, with the indifference price \eqref{european} for an European style claim used instead of a risk-neutral expectation at each node of the tree. Similarly, the employee exit rate can be naturally incorporated in the model by assuming that it is given by an independent Poisson process with intensity $\lambda$, so that the value of the position obtained in each node should be multiplied by the survival probability 
$(1-\lambda\Delta t)$ for each time step. We do not expect that either effect can change the nonlinearity and correlation features investigate in this paper.

\end{document}